\documentclass[12pt]{article}
\usepackage{amssymb,amsmath}
\usepackage{amsthm}
\usepackage{cases}
\usepackage{amsfonts}
\usepackage{graphicx}
\usepackage{makecell}
\usepackage{cite,color,xcolor}
\usepackage[margin=1.3 in]{geometry}
\usepackage[colorlinks,citecolor=blue,urlcolor=blue]{hyperref}
\usepackage[utf8]{inputenc}
\usepackage[pagewise]{lineno}

\newtheorem{theorem}{Theorem}[section]
\newtheorem{corollary}{Corollary}[section]
\newtheorem{lemma}{Lemma}[section]

\newtheorem{definition}{Definition}[section]
\newtheorem{remark}{Remark}[section]

\usepackage{txfonts}
\newcommand{\bal}{\begin{align}}
\newcommand{\bbal}{\begin{align*}}
\newcommand{\beq}{\begin{equation}}
\newcommand{\eeq}{\end{equation}}
\newcommand{\bca}{\begin{cases}}
\newcommand{\eca}{\end{cases}}
\def\div{\mathord{{\rm div}}}
\newcommand{\pa}{\partial}
\newcommand{\fr}{\frac}

\newcommand{\De}{\Delta}

\newcommand{\ep}{\varepsilon}
\newcommand{\dd}{\mathrm{d}}

\newcommand{\R}{\mathbb{R}}
\newcommand{\T}{\mathbb{T}}

\newcommand{\f}{\left}
\newcommand{\g}{\right}

\numberwithin{equation}{section}

\begin{document}

\title{Ill-posedness and inviscid limit of the basic equations of fluid dynamics in Besov spaces}

\author{
Yanghai Yu\footnote{
 School of Mathematics and Statistics,
 Anhui Normal University, Wuhu 241002, China.
\text{E-mail: yuyanghai214@sina.com}}
\qquad
 Jinlu Li\footnote{
 School of Mathematics and Computer Sciences,
 Gannan Normal University, Ganzhou 341000, China.
\text{E-mail: lijinlu@gnnu.edu.cn}}
\qquad
 Xing Wu\footnote{
College of Artificial Intelligence,
Henan Agricultural University, Zhengzhou 450046, China.
\text{E-mail: ny2008wx@163.com}}
}
\date{}
\maketitle
\begin{abstract}
In this paper, we consider the Cauchy problem to the basic equations of fluid dynamics on the torus. Firstly, we construct a new initial data and provide a simple proof on the ill-posedness of $B^s_{p,\infty}$ solution of the Euler equations and the surface quasi-geostrophic equations, which covers the results obtained by Cheskidov-Shvydkoy \cite{CS} and Misio{\l}ek-Yoneda \cite{MY}. Secondly, we prove the failure of the $B^s_{p,\infty}$-convergence in the inviscid limit for both the Navier-Stokes equations and the surface quasi-geostrophic equations.
\end{abstract}


\section{Introduction}
In this paper, we consider the Cauchy problem for the incompressible Navier-Stokes equations in $\T^d$ with $d\geq2$
\begin{align}\tag{NS}\label{ns}
\begin{cases}
\pa_t u+u\cdot \nabla u-\ep\Delta u+\nabla P=0, &\quad (t,x)\in \R^+\times\T^d,\\
\mathrm{div\,} u=0,&\quad (t,x)\in \R^+\times\T^d,\\
u(0,x)=u_0(x), &\quad x\in \T^d=(\R\setminus2\pi \mathbb{Z})^d,
\end{cases}
\end{align}
where $\ep>0$ is the kinematic viscosity, the vector field $u(t,x):[0,\infty)\times {\mathbb T}^d\to {\mathbb R}^d$ stands for the velocity of the fluid, the quantity $P(t,x):[0,\infty)\times {\mathbb T}^d\to {\mathbb R}$ denotes
the scalar pressure, and $\mathrm{div\,} u=0$ means that the fluid is incompressible. The Navier-Stokes systems describe various models of fluid
mechanics, engineering, biology, chemistry, and there is a huge list of references
concerning the mathematical analysis of related boundary value problems and
of their applications \cite{AA2}.

When the viscosity vanishes ($\ep=0$), the Navier-Stokes equations \eqref{ns} reduce to the Euler equations for ideal incompressible fluid
\begin{align}\tag{E}\label{e}
\begin{cases}
\pa_t u+u\cdot \nabla u+\nabla P=0, &\quad (t,x)\in \R^+\times\T^d,\\
\mathrm{div\,} u=0,&\quad (t,x)\in \R^+\times\T^d,\\
u(0,x)=u_0(x), &\quad x\in \T^d=(\R\setminus2\pi \mathbb{Z})^d.
\end{cases}
\end{align}
We say that the Cauchy problem \eqref{ns} or \eqref{e} is Hadamard (locally) well-posed in a Banach space $X$ if for any data $u_0\in X$ there exists (at least for a short time) $T>0$ and a unique solution in the space $\mathcal{C}([0,T),X)$ which depends continuously on the data. In particular, we say that the solution map is continuous if for any $u_0\in X$, there exists a neighborhood $B \subset X$ of $u_0$ such that
for every $u \in B$ the map $u \mapsto U$ from $B$ to $\mathcal{C}([0, T]; X)$ is continuous, where $U$
denotes the solution to \eqref{ns} or \eqref{e} with initial data $u_0$. For the well-posedness of \eqref{e} in Besov spaces we refer to \cite[Theorem 7.1]{BCD} and \cite[Theorem 1.1]{guo}.

Next, we review some developments which are closely related to the present work.
Kato \cite{Kato} obtained the local well-posedness of classical solution to Euler equations in the Sobolev space $H^s(\mathbb{R}^3)$ for all $s>5/2$. Kato-Ponce \cite{KatoP} extended this result to the Sobolev spaces $W^{s, p}(\mathbb{R}^3)$ of the fractional order for $s>3 / p+1,1<p<\infty$. Chae \cite{Chae1,Chae2,Chae3} and Chen-Miao-Zhang \cite{CMZ} gave further extensions to the Triebel-Lizorkin spaces $F_{p, r}^s(\mathbb{R}^3)$ with $s>3 / p+1,1<p, r<\infty$ and the Besov spaces $B_{p, r}^s(\mathbb{R}^3)$ with $s>3 / p+1$, $1<p<\infty, 1 \leq r \leq \infty$ or $s=3 / p+1,1<p<\infty, r=1$ (see \cite{Sawano} for some general Banach spaces). However, these two kinds of function spaces are only in the $L^p(1<p<\infty)$-framework since the Riesz transform is not bounded on $L^{\infty}$. The currently-known best result on the local existence was given by Pak-Park \cite{Pak} in the Besov space $B_{\infty, 1}^1(\mathbb{R}^3)$. Guo-Li-Yin \cite{guo} proved the continuous dependence of the Euler equations in the space $B_{p, r}^s(\mathbb{R}^3)$ with $s>3 / p+1$, $1\leq p\leq \infty, 1 \leq r < \infty$ or $s=3 / p+1,1\leq p\leq \infty, r=1$. Cheskidov-Shvydkoy \cite{CS} proved that the solution of the Euler equations cannot be continuous as a function
of the time variable at $t = 0$ in the spaces $B^s_{p,\infty}(\mathbb{T}^d)$ where $s > 0$ if $2 < p\leq \infty$ and $s>d(2/p-1)$
if $1 \leq p \leq 2$. For the whole space case, Li-Yu-Zhu \cite{Lyz1} constructed an initial data $u_0\in B^s_{p,\infty}$ showing that the corresponding solution of the Euler equations starting from $u_0$ is discontinuous at $t = 0$ in the metric of $B^s_{p,\infty}(\R^d)$ with $s\geq1+d/p$ (see also \cite{YuLf}). Bourgain-Li in \cite{B1,B2} employed a combination of Lagrangian and Eulerian techniques to obtain strong local ill-posedness results in borderline Besov spaces $B^{d/p+1}_{p,r}$ for $1\leq p<\infty$ and $1<r\leq\infty$ when $d=2,3$.
Misio{\l}ek-Yoneda \cite{MY2} showed that the solution for the Euler equations is not even continuous in the space of H\"{o}lder continuous functions and thus not locally Hadamard well-posed in $C^{1,\alpha}=B^{1+\alpha}_{\infty,\infty}$ with any $\alpha\in(0,1)$.

In this paper, we consider the ill-posedness problem of the Euler equations in Besov spaces.
The first result of this paper reads as follows:
\begin{theorem}\label{th1}
Let $d\geq2$ and $s>0$ with $1\leq p\leq \infty$. There exist initial data $u_0\in B^s_{p,\infty}(\T^d)$ and a corresponding weak solution $u$
of the Euler equations \eqref{e}  such that
\bbal
\limsup_{t\to0^+}\f\|u(t)-u_0\g\|_{B^s_{p,\infty}(\T^d)}\geq \eta_0,
\end{align*}
for some positive constant $\eta_0$ dependent on $p$ and $d$.
\end{theorem}
\begin{remark}\label{rem1}
Theorem \ref{th1} demonstrates the ill-posedness of the Euler equations in $B^s_{p,\infty}$. More precisely, there exist $u_0\in B^s_{p,\infty}$ and a corresponding weak solution to the Euler equations that starts from $u_0$ does not converge back to $u_0$ in the metric of $B^s_{p,\infty}$ as time goes to zero.
\end{remark}
\begin{remark}\label{rem2}
Theorem \ref{th1} improves the result obtained by Cheskidov and Shvydkoy in \cite{CS} since the index $p>2$ has been enlarged to $1\leq p\leq \infty$.
\end{remark}
\begin{remark}\label{rem3}
Compared with the idea in Cheskidov and Shvydkoy in \cite{CS}, we present a short and direct proof by constructing explicit traveling wave solutions, which does not depend on the nonlinear structure of the Euler equations. In fact, we construct special solutions causing the ill-posedness. To clearly see this, for instance, let $s>0$ and $n$ be a sufficiently large number, we can
verify that
\begin{align*}
v(x_1-t)=2^{-ns}\cos \f(\frac{11}{8}2^{n}(x_1-t)\g)
\quad\text{and}\quad
v(x_1)=2^{-ns}\cos \f(\frac{11}{8}2^{n}x_1\g)
\end{align*}
are two high frequency wave. In this case, both $\left\|v(x_1-t)\right\|_{B_{p, \infty}^s(\mathbb{T}^d)}$ and $\left\|v(x_1)\right\|_{B_{p, \infty}^s(\mathbb{T}^d)}$ have  positive lower and supper bounds. However, we notice that
$v(x_1-t)=-v(x_1)$ if $\frac{11}{8}2^{n}t=\pi$,
thus the Besov norm of the difference, i.e., $$\left\|v(x_1-t)-v(x_1)\right\|_{B_{p, \infty}^s(\mathbb{T}^d)}=2\left\|v(x_1)\right\|_{B_{p, \infty}^s(\mathbb{T}^d)},$$ still has a positive lower bound.
\end{remark}

We recall the classical well-posedness theory for the Euler equations, as founded in \cite{BCD}, which is convenient for our purposes: Given a divergence-free initial data $u_0 \in B_{p, r}^s$ with $s>1+d / p,1 \leq p, r \leq \infty$ or $s=1+d / p,1 \leq p \leq \infty,r=1$, then there exists a short time $T=T(u_0)$, such that \eqref{e} has a unique solution $u(t,x)\in C([0, T] ; B_{p, r}^s)$ if $1\leq r<\infty$ or $u(t,x)\in C_w([0, T] ; B_{p, \infty}^s)$, where the subscript $w$ indicates weak continuity in the time variable. Furthermore, we have
$\|u(t)\|_{B_{p, r}^s} \lesssim\left\|u_0\right\|_{B_{p, r}^s}$ for $0 <t \leq T$.

From now on, we denote
 \begin{align*}
&u^{\rm E}(t,u_0)=\text{\rm the solution map of the Euler equations \eqref{e}}\, \text{\rm with initial data}\, u_0
 \end{align*}
and denote any bounded subset $U_R$ in $B^s_{p,r}(\mathbb{T}^d)$ by
$$U_R:=\left\{u_0\in B^s_{p,r}(\mathbb{T}^d): \|u_0\|_{B^s_{p,r}(\mathbb{T}^d)}\leq R,\;\mathrm{div}\, u_0=0\right\}.$$

Our second result of this paper reads as follows:
 \begin{theorem}\label{th2}
Let $d\geq2$. Assume that $(s,p,r)$ satisfies
\begin{align}\label{spr}
s>\frac{d}{p}+1, (p,r)\in [1,\infty]\times [1,\infty) \quad   \mathrm{or}    \quad s=\frac{d}{p}+1, (p,r)\in [1,\infty]\times\{1\}.
\end{align}
For any $\alpha\in(0,1)$, there exists initial data $u_0\in U_R$ such that the data-to-solution map $u_0\mapsto u^{\rm E}(t,u_0)\in B^s_{p,r}(\T^d)$
of the Euler equations \eqref{e} satisfies
\bbal
\limsup_{t\to0^+}\frac{\f\|u^{\rm E}(t,u_0)-u_0\g\|_{B^s_{p,r}(\T^d)}}{t^\alpha}=+\infty.
\end{align*}
\end{theorem}

\begin{remark}\label{r3}
It is known that the solution $u^{\rm E}(t,u_0)$ for \eqref{e} is continuous in time
in Besov spaces $B^s_{p,r}$ with $(s,p,r)$ satisfies \eqref{spr}, while Theorem \ref{th2} furthermore indicates that the solution $u^{\rm E}(t,u_0)$ for \eqref{e} cannot be H\"{o}lder continuous in time in the same Besov spaces.
\end{remark}
A classical problem in fluid mechanics is the approximation in the limit $\ep\to0$ of vanishing viscosity (also called inviscid limit) of solutions of the Euler equations by solutions of the incompressible Navier-Stokes equations.
The problem of the convergence of smooth viscous solutions of \eqref{ns} to the Eulerian one as $\ep\to0$ is
well understood and has been studied in many literature, see for example \cite{C17,Swann, Kato}, and \cite{CKV,CV} for the inviscid limit in a bounded domain. Majda \cite{Majda} showed that under the assumption $u_0\in H^s$ with $s>\frac{d}{2}+2$, the solutions $u_\ep$ to \eqref{ns} converge in $L^2$ norm as $\ep\to 0$ to the unique solution of Euler equations and the convergence rate is of order $(\ep t)^{\fr12}$. Masmoudi \cite{M} improved the result and obtained the convergence in $H^s$-norm under the assumption $u_0\in H^s$ with $s>\frac{d}{2}+1$. In dimension two, Hmidi and Kerrani in \cite{HK} proved that \eqref{ns} is globally well-posed in the Besov space $B^2_{2,1}$, with uniform bounds on the viscosity and obtained that the convergence rate of the inviscid limit is of order $\ep t$ for
vanishing viscosity. Subsequently, in \cite{HK1}, they further generalized to other Besov spaces $B^{2/p+1}_{p,1}$ with convergence in $L^p$. Chemin \cite{chem} resolved inviscid limit of Yudovich type solutions with only the assumption that the vorticity is bounded. In the case of $\T^2$ or $\R^2$, by taking greater advantage of vorticity formulation, more beautiful results were obtained quantitatively (see for example \cite{Be,Ci,P2} and the references therein). Guo-Li-Yin \cite{guo} solved the inviscid limit in the same topology. However, it left an open problem for the end-point case $r=\infty$ in the inviscid limit of the Navier-Stokes equations in Besov spaces. For the whole space case, Yu and Li \cite{YuL} proved the failure of the $B^{s}_{p,\infty}$-convergence of the Navier-Stokes equations toward the Euler equations in the inviscid limit. In this paper, we consider the torus case and establish the following result

\begin{theorem} {\bf (Non-convergence)} \label{th4} Let $d\geq 2$ and $\ep\in [0,1]$. Assume that $s>1+d/p$ with $1\leq p\leq \infty$. Then a family of solution maps $(u_0,\ep)\mapsto u^{\rm NS}_{\ep}(t,u_0)$ to the Navier-Stokes equations \eqref{ns}
do not converge to the solution map $u_0\mapsto u^{\rm E}(t,u_0)$ of the Euler equations \eqref{e} in $B^{s}_{p,\infty}$.
More precisely, there exists initial data $u_0\in U_R$ such that
$$
\limsup_{\ep_n\to 0^+}\left\|u^{\rm NS}_{\ep_n}(t_n,u_0)-u^{\rm E}(t_n,u_0)\right\|_{B^s_{p,\infty}(\T^d)}\geq \eta_0,
$$
for some positive constant $\eta_0$ only dependent on $p$ and $d$.
\end{theorem}
\begin{remark}
We would like to mention that Theorems \ref{th1}-\ref{th4} are new. In fact, compared with the non-periodic results in \cite{Lyz1,YuLf,BULL,YuL}, we consider the torus case and provide a simple proof of Theorems \ref{th1}-\ref{th4} in a unified way by constructing explicit periodic solutions. In particular, Theorem \ref{th1} improves the regularity index $s\geq1+d/p$ in \cite{Lyz1,YuLf} to $s>0$.

\end{remark}
Next, we consider the Cauchy problem for the two-dimensional surface quasi-geostrophic equations which is a fundamental example of active scalar transport
\begin{align}\tag{QG}\label{QG}
\begin{cases}
\pa_t \theta+u\cdot \nabla \theta+\kappa\Lambda^{\alpha}\theta=0, &\quad (t,x)\in \R^+\times\T^2,\\
u=\mathcal{R}^{\perp} \theta:=\nabla^{\perp} \Lambda^{-1} \theta,&\quad (t,x)\in \R^+\times\T^2,\\
\theta(0,x)=\theta_0(x), &\quad x\in \T^2=(\R\setminus2\pi \mathbb{Z})^2,
\end{cases}
\end{align}
where the unknown $\theta=\theta(x, t)$ is scalar field,
$\Lambda^\alpha f$ is defined via the Fourier modes of $f$
\begin{align*}
\widehat{\Lambda^\alpha f}(k)=|k|^{\alpha} \widehat{f}(k).
\end{align*}
$\mathcal{R}=\left(\mathcal{R}_1, \mathcal{R}_2\right)$ is the vector of Riesz transforms, then
$$u=\left(\mathcal{R}_2 \theta,-\mathcal{R}_1 \theta\right), \quad \widehat{\mathcal{R}_j \theta}(\xi)=-\mathrm{i} \frac{\xi_j}{|\xi|} \hat{\theta}(\xi), \quad \xi_j \in \mathbb{Z}, \quad j=1,2.$$
When the viscosity vanishes ($\kappa=0$), \eqref{QG} reduces to the inviscid quasi-geostrophic equations. The QG equations were introduced in \cite{CMT94} as a two-dimensional model of the 3D Euler equations. There have been significant development in the ill-posedness theory (see e.g.,\cite{CM22ADV,CM24CMP,JK24APDE}) for the QG equations. To put our study in the proper perspective, we first recall a few results which are related to our problem.

\underline{Ill-posedness in $W^{1, \infty}$}: Elgindi and Masmoudi \cite{Elgindi20arma} proved the ill-posedness of active scalar system in $L^{\infty}$-type spaces in the sense that there exist smooth steady states $\bar{\theta}$ and a sequence of perturbations $\tilde{\theta}_0^{(\epsilon)}\left(\epsilon \rightarrow 0^{+}\right)$ such that the associated solution $\theta^{(\epsilon)}$ with data $\bar{\theta}+\tilde{\theta}_0^{(\epsilon)}$ satisfies
$$
\left\|\theta^{(\epsilon)}(0, \cdot)-\bar{\theta}\right\|_{W^{1, \infty}}<\epsilon, \quad \sup _{0<t<\epsilon}\left\|\theta^{(\epsilon)}(t, \cdot)-\bar{\theta}\right\|_{W^{1, \infty}}>c(\bar{\theta}).
$$

\underline{Ill-posedness in $B_{p, \infty}^s$}: Kiselev, Nazarov and Volberg \cite{KNV} observed that if $\theta_0 \in L^p(\mathbb{T}^2)$ with $1<p<\infty$, the solution of \eqref{QG} with $\alpha=1$ satisfies $\lim _{t \rightarrow 0^{+}}\left\|\theta(t)-\theta_0\right\|_{L^p}=0$. Misio{\l}ek and Yoneda \cite{MY} showed that this property fails in certain Besov spaces for $0<\alpha<1$ by using elementary properties of certain lacunary Fourier series.
\begin{theorem}[\cite{MY}] Let $\alpha\in (0,1)$. There exist initial data $\theta_0 \in B_{p, \infty}^s(\mathbb{T}^2)$ and a corresponding weak solution $\theta$ of the Cauchy problem \eqref{QG} with $\kappa>0$ such that
$$
\lim _{t \rightarrow 0^{+}}\left\|\theta(t)-\theta_0\right\|_{B_{p, \infty}^s(\T^2)}>0
$$
for any $s>0$ and any $2 \leq p \leq \infty$.
\end{theorem}
 Our aim is to generalize the above result to the broader cases.

\begin{theorem}[$\kappa>0$]\label{Qth1}
Let $\alpha\in(0,2]$, $s>0$ and $1\leq p\leq \infty$. There exist initial data $\theta_0\in B^s_{p,\infty}(\T^2)$ and a corresponding weak solution $\theta$
of the Cauchy problem \eqref{QG} with $\kappa>0$ such that
\bbal
\limsup_{t\to0^+}\f\|\theta(t)-\theta_0\g\|_{B^{s}_{p,\infty}(\T^2)}\geq \eta_0,
\end{align*}
for some positive constant $\eta_0$ which is only dependent on $p$.
\end{theorem}
\begin{remark}\label{r5}
Theorem \ref{Qth1} demonstrates the ill-posedness of the QG equations in $B^s_{p,\infty}$ and covers the result obtained by Misio{\l}ek and Yoneda \cite{MY}. More precisely, there exist $u_0\in B^s_{p,\infty}$ and a corresponding weak solution to the QG equations that starts from $u_0$ does not converge back to $u_0$ in the metric of $B^s_{p,\infty}$ as time goes to zero. We should mention that
Theorem \ref{Qth1} holds for the generalised QG equations (see e.g.,\cite{CM22ADV,CM24CMP})
\begin{align*}
\begin{cases}
\pa_t \theta+u\cdot \nabla \theta+\kappa\Lambda^{\alpha}\theta=0, \\
u=\nabla^{\bot}\Lambda^{\beta-2}\theta,\\
\theta(0,x)=\theta_0(x).
\end{cases}
\end{align*}
\end{remark}
\begin{theorem}\label{Qth2} {\bf (Non-convergence)} Let $\alpha\in(0,2]$ and $s>1+\fr{2}{p}$ with $1\leq p\leq \infty$.
There exists initial data $\theta_0\in B^{s}_{p,\infty}(\mathbb{T}^2)$ such that
the solution $\theta^{\kappa}(t)$ of the dissipative \eqref{QG} does not converge to the solution $\theta^{0}(t)$ of the inviscid \eqref{QG} for small $t\in(0,T_1]$ in $B^{s}_{p,\infty}$ as $\kappa\downarrow 0$.
More precisely, there exists initial data $\theta_0\in B^{s}_{p,\infty}(\mathbb{T}^2)$ such that for a short time $t=\kappa$
\begin{align*}
\limsup_{\kappa\to 0^+}\left\|\theta^{\kappa}(t,\theta_0)-\theta^{0}(t,\theta_0)\right\|_{B^{s}_{p,\infty}(\mathbb{T}^2)}\geq \eta_0,
\end{align*}
for some positive constant $\eta_0$ depending on $p$ but independent on $\kappa$.
\end{theorem}

\begin{theorem}[$\kappa=0$]\label{Qth3}
Let $s>0$ and $1\leq p\leq \infty$. There exist initial data $\theta_0\in B^s_{p,\infty}(\T^2)$ and a corresponding weak solution $\theta$
of the inviscid \eqref{QG} such that for some positive constant $\eta_0$ which is dependent on $p$
\bal\label{imp2}
\limsup_{t\to0^+}\f\|\theta(t)-\theta_0\g\|_{B^{s}_{p,\infty}\cap L_{\mathcal{R}}^{\infty}(\T^2)}\geq \eta_0,
\end{align}
where
\bbal
\f\|f\g\|_{B^{s}_{p,\infty}\cap L_{\mathcal{R}}^{\infty}}:=\f\|f\g\|_{B^{s}_{p,\infty}}+\f\|f\g\|_{L_{\mathcal{R}}^{\infty}}\quad\text{with}\quad
\f\|f\g\|_{L_{\mathcal{R}}^{\infty}}:=\f\|f\g\|_{L^{\infty}}+\|\mathcal{R}^{\perp}f\|_{L^{\infty}}.
\end{align*}
\end{theorem}

\begin{corollary}[$\kappa=0$]\label{c2}
Let $s>\fr{2}{p}$ and $1\leq p< \infty$. There exist initial data $\theta_0\in B^s_{p,\infty}(\T^2)$ and a corresponding weak solution $\theta$ of the inviscid \eqref{QG} such that
\bbal
\limsup_{t\to0^+}\f\|\theta(t)-\theta_0\g\|_{B^{s}_{p,\infty}(\T^2)}\geq \eta_0,
\end{align*}
for some positive constant $\eta_0$ which is dependent on $p$.
\end{corollary}
\begin{remark}\label{r18}
Theorem \ref{Qth3} implies the ill-posedness of the inviscid \eqref{QG} in $B^{s}_{p,\infty}\cap L_{\mathcal{R}}^{\infty}$ with $s>0$
while Corollary \ref{c2} implies the ill-posedness of the inviscid \eqref{QG} in $B^s_{p,\infty}$ with $s>\fr{2}{p}$.
The comparison between the proof of Theorem \ref{Qth1} and Theorem \ref{Qth3} tells us that the mechanisms underlying the discontinuity of solution in weaker Besov spaces are completely different for the dissipative and inviscid QG equations. Precisely speaking, the primarily affect which leads to the ill-posedness of the dissipative QG equations is the diffusion term while for the inviscid QG equations it is the convection term.
\end{remark}

\section{Preliminaries}\label{sec2}
We define the periodic Fourier transform $\mathcal{F}_{\mathbb{T}^d}: \mathcal{D}(\mathbb{T}^d)\rightarrow \mathcal{S}(\mathbb{Z}^d)$  as
$$(\mathcal{F} u)(k)=\widehat{u}(k)=\frac{1}{(2 \pi)^d} \int_{\mathbb{T}^d} e^{-\mathrm{i} x \cdot k} u(x) \mathrm{d} x.$$
We decompose $u \in \mathcal{D}(\mathbb{T}^d)$ on the circle $\mathbb{T}^d$ into Fourier series, i.e.
$$
u(x)=\sum_{k \in \mathbb{Z}^d} \widehat{u}(k) e^{\mathrm{i} x \cdot k}.
$$
We are interested in solutions which take values in the Besov space $B_{p, r}^s(\mathbb{T}^d)$. Recall that one way to define this space requires a dyadic partition of unity. Given a smooth bump function $\chi$ supported on the ball of radius $4 / 3$, and equal to 1 on the ball of radius $3 / 4$, we set $\varphi(\xi)=\chi(2^{-1} \xi)-\chi(\xi)$ and $\varphi_j(\xi)=\varphi(2^{-j} \xi)$ and then we deduce that $\varphi$ satisfies that ${\rm{supp}} \;\varphi\subset \left\{\xi\in \T^d: 3/4\leq|\xi|\leq8/3\right\}$ and
 $\varphi(\xi)\equiv 1$ for $4/3\leq |\xi|\leq 3/2$. Using this partition, we define the periodic dyadic blocks as follows
\begin{align*}
&\Delta_j u=0, \, \text{ if } \; j \leq-2, \\
&\Delta_{-1} u=\sum_{\xi \in \mathbb{Z}^d} \chi(\xi) \widehat{u}(\xi) e^{\mathrm{i} x \cdot \xi}, \\
&\Delta_j u= \sum_{\xi\in \mathbb{Z}^d} \varphi_j( \xi)\widehat{u}(\xi) e^{\mathrm{i} x \cdot \xi}, \, \text {if} \; j \geq 0 .
\end{align*}
The operators $\Delta_j$ are defined on the periodic domain and shares many properties with those
on the whole space (see \cite{BCD,Dai20}). In particular, we obtain the Littlewood-Paley decomposition of $u$
$$u=\sum_{j\geq-1} \Delta_j u\quad \text{in}\quad \mathcal{S}^{\prime}(\mathbb{T}^d).$$

\begin{definition}
Let $s\in\mathbb{R}$ and $(p,r)\in[1, \infty]^2$. The nonhomogeneous Besov space $B^{s}_{p,r}(\mathbb{T}^d)$ is defined by
\begin{align*}
B^{s}_{p,r}(\mathbb{T}^d):=\f\{f\in \mathcal{S}'(\mathbb{T}^d):\;\|f\|_{B^{s}_{p,r}(\mathbb{T}^d)}<\infty\g\},
\end{align*}
where
\begin{numcases}{\|f\|_{B^{s}_{p,r}(\mathbb{T}^d)}=}
\left(\sum_{j\geq-1}\left(2^{sj}\|\Delta_jf\|_{L^p(\mathbb{T}^d)}\right)^r\right)^{1/r}, &\text {if} \; $1\leq r<\infty$,\nonumber\\
\sup_{j\geq-1}\left(2^{sj}\|\Delta_jf\|_{L^p(\mathbb{T}^d)}\right), &\text {if} \;  $r=\infty$.\nonumber
\end{numcases}
\end{definition}

Let us complete this section by presenting two lemmas which will be used often in the sequel.
\begin{lemma}\label{ley1}
Let $3\leq m\in \mathbb{Z}$ and $-1\leq j\in \mathbb{Z}$, we have
$$\Delta_j\f[\cos \f(\frac{11}{8}2^{m} x_1\g)\g]=\begin{cases}
0, &\text{if}\; j\neq m,\\
\cos\f(\frac{11}{8}2^{m} x_1\g), &\text{if}\; j=m.\end{cases}$$
\end{lemma}
\begin{proof}\; We set $\lambda=\frac{11}{8}2^{m}$ and $\vec{\lambda}=(\lambda,0,\cdots,0)$ for simplicity. Notice that $\varphi_j(k)$ is symmetric, i.e., $\varphi_j(k)=\varphi_j(|k|)$, we deduce
\begin{align*}
\Delta_j\f(\cos (\lambda x_1)\g)&=\fr12\sum_{k\in \mathbb{Z}^d}\varphi_j(k)\f(\mathbf{1}_{\vec{\lambda}}(k)+\mathbf{1}_{-\vec{\lambda}}(k)\g)e^{\mathrm{i} x\cdot k}
\\&=\fr12\sum_{\ell\in \mathbb{Z}}\varphi_j(\ell)\f(\mathbf{1}_{\lambda}(\ell)+\mathbf{1}_{-\lambda}(\ell)\g)e^{\mathrm{i} x_1\ell}\\
&=\varphi_j(\lambda)\cos (\lambda x_1)=\begin{cases}
0, &\text{if}\; j\neq m,\\
\cos (\lambda x_1), &\text{if}\; j=m,
\end{cases}
\end{align*}
where $\mathbf{1}_{K}(x)$ is the indicator function, taking a value of 1 if $x=K$ and 0 otherwise.

Thus we obtain the desired result of Lemma \ref{ley1}.
\end{proof}
\begin{remark}\label{re6}
For large fixed $m$, $\varphi_j\f(\frac{11}{8}2^{m} x_1\g)$ vanish except for $j=m$ due to the support condition of $\varphi$ and $\varphi(k)\equiv 1$ for $4/3\leq |k|\leq 3/2$.
\end{remark}
The following simple fact is needed in the sequel.
\begin{lemma}\label{ley2} Let $1\ll n\in \mathbb{Z}$, we have
\bbal
\left\|\cos\left(\frac{11}{8}2^{n}x\right)\right\|_{L^p([0,2\pi])}=c_0:=\begin{cases}
\left(2\int^\pi_0|\cos x|^p\dd x\right)^{1/{p}}, &\text{if}\quad p\in[1,\infty),\\
1, &\text{if}\quad p=\infty.
\end{cases}
\end{align*}
\end{lemma}
\section{Proof of Theorem \ref{th1}}\label{sec3}
\noindent{\bf Choice of initial data.}\;Let us fix an $s>0$ and define the initial data
\bal\label{i}
u_{0}(x)=(1,f(x_1),0,\ldots,0),
\end{align}
where $f(x)$ is a bounded real-valued periodic function of one variable with the following form
\bal\label{f}
f(x)=\sum\limits^\infty_{j=3}2^{-js}\cos \f(\frac{11}{8}2^{j}x\g).
\end{align}
It is not difficult to check that $\div\, u_0=0$ and $u_0\in B^s_{p,\infty}(\T^d)$ for any $1\leq p\leq \infty$.

\noindent{\bf Construction of solution.}\;An fundamental observe is  that the vector field
\bal\label{s1}
&u(t,x)=\f(1,f\f(x_1-t\g),0,\ldots,0\g)
\end{align}
is an obvious periodic solution of the incompressible Euler equations:
\bbal\pa_t u+u\cdot \nabla u+\nabla P=0, \quad
\mathrm{div\,} u=0,
\end{align*}
with $P = 0$, i.e., this is a pressureless flow.

From \eqref{i} and \eqref{s1}, we have
\bbal
&u(t,x)-u_{0}(x)=\f(0,f\f(x_1-t\g)-f(x_1),0,\ldots,0\g),
\end{align*}
From the explicit formula for $f(x)$ in \eqref{f} and by direct computation, we see that
\bbal
f(x_1-t)-f(x_1)&=\sum\limits^\infty_{j=3}2^{-js}\f[\cos \f(\frac{11}{8}2^{j}x_1\g)\mathbf{a}_j(t)+\sin \f(\frac{11}{8}2^{j}x_1\g)\mathbf{b}_j(t)\g],
\end{align*}
where
\bal\label{ab}
&\mathbf{a}_j(t)=\cos \f(\frac{11}{8}2^{j}t\g)-1\quad \text{and}\quad \mathbf{b}_j(t)=\sin \f(\frac{11}{8}2^{j}t\g).
\end{align}
Then, by Lemma \ref{ley1}, we have for some $n$ large enough
\bbal
2^{ns}\De_n[f(x_1-t)-f(x_1)]&=\cos \f(\frac{11}{8}2^{n}x_1\g)\mathbf{a}_n(t)+\sin \f(\frac{11}{8}2^{n}x_1\g)\mathbf{b}_n(t).
\end{align*}
Letting $\frac{11}{8}2^{n}t_n=\pi$, then one has $\mathbf{a}_n(t_n)=-2$ and $\mathbf{b}_n(t_n)=0$.
Thus we have for some $n$ large enough
\bbal
&\quad 2^{ns}\f\|\De_n[f(x_1-t_n)-f(x_1)]\g\|_{L^p(\T^d)}
=2(2\pi)^{\frac{d-1}{p}}\f\|\cos \f(\frac{11}{8}2^{n}x_1\g)\g\|_{L^p(\mathbb{T})}.
\end{align*}
From which and Lemma \ref{ley2}, we deduce that for some $n$ large enough
\bbal
\|u(t_n,x)-u_{0}(x)\|_{B^s_{p,\infty}(\T^d)}&=\|f(x_1-t_n)-f(x_1)\|_{B^s_{p,\infty}(\T^d)}\\
&\geq2^{ns}\f\|\De_n[f(x_1-t_n)-f(x_1)]\g\|_{L^p(\T^d)}\\&= 2c_0(2\pi)^{\frac{d-1}{p}}.
\end{align*}
Notice that $t_n\to 0^+$ as $n\to\infty$, we complete the proof of Theorem \ref{th1}.

\section{Proof of Theorem \ref{th2}}\label{sec4}

Assume that $(s,p,r)$ satisfies \eqref{spr}.
Let us define the initial data
$$u_{0}(x)=(1,g(x_1),0,\cdots,0),$$
where $g(x)$ is a bounded real-valued period function of one variable with the following form
\bal\label{g}
g(x)=\sum\limits^\infty_{j=3}j^{-2}2^{-js}\cos \f(\frac{11}{8}2^{j}x\g).
\end{align}
It is not difficult to check that $\div\, u_0=0$ and $u_0\in B^s_{p,r}(\T^d)$.

Let $t\geq 0$, we consider
\bbal
&u(t,x)=\f(1,g\f(x_1-t\g),0,\ldots,0\g),
\end{align*}
Obviously, it is a period function which also satisfies the incompressible Euler equations \eqref{e} with initial condition $u_{0}(x)=(1,g(x_1),0,\cdots,0)$.

Continue in a similar fashion we see that
\bbal
&u(t,x)-u_{0}(x)=\f(0,g\f(x_1-t\g)-g(x_1),0,\ldots,0\g),
\end{align*}
and
\bbal
g(x_1-t)-g(x_1)&=\sum\limits^\infty_{j=3}j^{-2}2^{-js}\f[\cos \f(\frac{11}{8}2^{j}x_1\g)\mathbf{a_j}(t)+\sin \f(\frac{11}{8}2^{j}x_1\g)\mathbf{b}_j(t)\g],
\end{align*}
where $\mathbf{a_j}(t)$ and $\mathbf{b_j}(t)$ are given by \eqref{ab}.

Then, by Lemma \ref{ley1}, we have for some $n$ large enough
\bbal
n^{2}2^{ns}\De_n[g(x_1-t)-g(x_1)]&=\cos \f(\frac{11}{8}2^{n}x_1\g)\mathbf{a}_n(t)+\sin \f(\frac{11}{8}2^{n}x_1\g)\mathbf{b}_n(t).
\end{align*}
Letting $\frac{11}{8}2^{n}t_n=\pi$, then one has $\mathbf{a}_n(t_n)=-2$ and $\mathbf{b}_n(t_n)=0$. Thus we have for some $n$ large enough
\bbal
&\quad 2^{ns}\f\|\De_n[g(x_1-t_n)-g(x_1)]\g\|_{L^p(\T^d)}
=2(2\pi)^{\frac{d-1}{p}}n^{-2}\f\|\cos \f(\frac{11}{8}2^{n}x_1\g)\g\|_{L^p(\mathbb{T})}.
\end{align*}
From which and Lemma \ref{ley2}, we deduce that for some $n$ large enough
\bbal
t_n^{-\alpha}\|u(t_n,x)-u_{0}(x)\|_{B^s_{p,r}(\T^d)}&= t_n^{-\alpha}\|g(x_1-t_n)-g(x_1)\|_{B^s_{p,r}(\T^d)}\\
&\geq t_n^{-\alpha}2^{ns}\f\|\De_n[g(x_1-t_n)-g(x_1)]\g\|_{L^p(\T^d)}\\
&= 2c_0(2\pi)^{\frac{d-1}{p}}t_n^{-\alpha}n^{-2}.
\end{align*}
Notice that $t_n\to 0^+$ and $t_n^{-\alpha}n^{-2}\to+\infty$ as $n\to\infty$, we complete the proof of Theorem \ref{th2}.

\section{Proof of Theorem \ref{th4}}\label{sec6}
\noindent{\bf Choice of initial data.}\;Let us fix an $s$ satisfying $s>1+d/p$ with $1\leq p\leq \infty$ and define the initial data
\bbal
u_{0}(x)=(1,f(x_1),0,\ldots,0),
\end{align*}
where $f(x)$ is given by \eqref{f}.

\noindent{\bf Construction of solution to the Euler equations.}\;

Let $t\geq0$, we  consider
\bbal
&u^{\rm E}(t,x)=\f(1,f(x_1-t),0,\ldots,0\g).
\end{align*}
By the classical well-posedness theory, we know that $u^{\rm E}(t,x)$ is a unique periodic solution of the Euler equations \eqref{e} with initial data $u_0(x)=\f(1,f(x_1),0,\ldots,0\g)$.

\noindent{\bf Construction of solution to the Navier-Stokes equations.}\;

Assume that $f_n(t,x_1)$ solves the Cauchy Problem
\begin{align}\label{c}
\begin{cases}
 \pa_tf_n+\pa_{x_1}f_n-\ep_n\pa^2_{x_1}f_n=0,\quad t>0,\\
f_{n}(t=0,x_1)=f(x_1).
\end{cases}
\end{align}
Let $t\geq 0$, we consider
\bbal
&u^{\rm NS}_{\ep_n}(t,x)=(1,f_n(t,x_1),0,\ldots,0).
\end{align*}
Due to \eqref{c}, we know that  $u^{\rm NS}_n(t,x)$ is a unique periodic solution of the Navier-Stokes equations with initial data $u_0(x)=\f(1,f(x_1),0,\ldots,0\g)$.

It is easy to deduce that the Cauchy Problem \eqref{c} has a unique explicit solution
\bbal
f_n(t,x_1)=\sum\limits^\infty_{j=3}2^{-js}\f[\cos \f(\frac{11}{8}2^{j}(x_1-t)\g)e^{-\ep_n\frac{121}{64}2^{2j}t}\g].
\end{align*}
We shall compare the solution of Euler equations with that of Navier-Stokes equations. Obviously, $$u^{\rm NS}_{\ep_n}(t_n,u_0)-u^{\rm E}(t_n,u_0)=(0,f_n(t,x_1)-f(x_1-t_n),0,\ldots,0),$$
and
\bbal
f_n(t,x_1)-f(x_1-t)&=\sum\limits^\infty_{j=3}2^{-js}\cos \f(\frac{11}{8}2^{j}(x_1-t)\g)\f(e^{-\ep_n\frac{121}{64}2^{2j}t}-1\g).
\end{align*}
By Lemma \ref{ley1}, one has for some $n$ large enough
\bbal
\De_n\f[f_n(t,x_1)-f(x_1-t_n)\g]
=2^{-ns}\f[\cos \f(\frac{11}{8}2^{n}x_1\g)\mathbf{e}_n(t)
+\sin \f(\frac{11}{8}2^{n}x_1\g)\mathbf{f}_n(t)\g]\mathbf{g}_n(t),
\end{align*}
where
\bbal
&\mathbf{e}_n(t)=\cos \f(\frac{11}{8}2^{n}t\g),\quad
\mathbf{f}_n(t)=\sin \f(\frac{11}{8}2^{n}t\g)\quad\text{and}\quad
\mathbf{g}_n(t)=e^{-\ep_n\frac{121}{64}2^{2n}t}-1.
\end{align*}
Letting $\frac{11}{8}2^{n}t_n=\pi$ and $\ep_n=\frac{8}{11\pi} 2^{-n}$, which gives that
$$\f(\mathbf{e}_n(t_n),\;\mathbf{f}_n(t_n),\;\mathbf{g}_n(t_n)\g)=(-1,\;0,\;e^{-1}-1)$$ and thus
\bbal
\De_n\f[f_n(t_n,x_1)-f(x_1-t_n)\g]&=(1-e^{-1})2^{-ns}\cos \f(\frac{11}{8}2^{n}x_1\g).
\end{align*}
Combining the above, we have
\bbal
\|u^{\rm NS}_{\ep_n}(t_n,x)-u^{\rm E}(t_n,u_0)\|_{B^s_{p,\infty}(\T^d)}&\geq2^{ns}\f\|\De_n\f[f_n(t_n,x_1)-f(x_1-t_n)\g]\g\|_{L^p(\T^d)}
\\&=\f(1-e^{-1}\g)
\cdot\f\|\cos \f(\frac{11}{8}2^{n}x_1\g)\g\|_{L^p(\T^d)}\\
&=c_0\f(1-e^{-1}\g)(2\pi)^{\frac{d-1}{p}}.
\end{align*}
Taking the $\limsup$, we deduce that
$$
\limsup_{\ep_n\to 0^+}\left\|u^{\rm NS}_{\ep_n}(t_n,u_0)-u^{\rm E}(t_n,u_0)\right\|_{B^s_{p,\infty}}\geq c_0\f(1-e^{-1}\g)(2\pi)^{\frac{d-1}{p}}.
$$
Notice that $t_n\to 0^+$ and $\ep_n\to0^+$ as $n\to\infty$, we complete the proof of Theorem \ref{th4}.

\section{Proof of Theorems \ref{Qth1}-\ref{Qth2}}\label{sec7}
\noindent{\bf Choice of initial data.}\;Let us fix an $s>0$ and define the initial data which is a bounded real-valued periodic function of one variable with the following form
\bbal
\theta_{0}(\mathbf{x})=f(x_1)=\sum\limits^\infty_{j=3}2^{-js}\cos \f(\frac{11}{8}2^{j}x_1\g), \quad  \mathbf{x}=(x_1,x_2)\in \mathbb{T}^2.
\end{align*}
It is not difficult to check that $\theta_{0}\in B^s_{p,\infty}(\T^2)$ for any $1\leq p\leq \infty$.

\noindent{\bf Construction of solution.}\;We are interested in solutions of the form
\bbal
&\theta(t,\mathbf{x})=f\f(t,x_1\g),
\end{align*}
then it must have
$$u\cdot\nabla \theta=u_1\partial_1\theta+u_2\partial_2\theta=0.$$
In fact,
$$\partial_2\theta=0\quad\text{and}\quad u_1=-\Lambda^{-1}\partial_2\theta=0.$$
Thus
$
\theta(t,\mathbf{x})=e^{-\kappa \Lambda^{\alpha}t}\theta_0(x_1)
$
is a periodic solution of \eqref{QG}. We shall compare the solution of \eqref{QG} with the initial data. Obviously,
\bbal
\theta(t,\mathbf{x})-\theta_{0}(\mathbf{x})&=\f(e^{-\kappa \Lambda^{\alpha}t}-1\g)\theta_{0}(\mathbf{x}).
\end{align*}
By Lemma \ref{ley1}, one has for some $n$ large enough
\begin{align*}
\Delta_j\f[\f(e^{-\kappa \Lambda^{\alpha}t}-1\g)\cos \f(\frac{11}{8}2^{n} x_1\g)\g]
&=\begin{cases}
0, &\text{if}\; j\neq n,\\
\f(e^{-\kappa \f(\frac{11}{8}2^{n}\g)^{\alpha}t}-1\g)\cos\f(\frac{11}{8}2^{n} x_1\g), &\text{if}\; j=n,
\end{cases}
\end{align*}
which gives that
\bal\label{hhy}
&\De_n\f[\theta(t,\mathbf{x})-\theta_{0}(\mathbf{x})\g]=2^{-ns}\cos \f(\frac{11}{8}2^{n}x_1\g)\f(e^{-\kappa\f(\frac{11}{8}2^{n}\g)^\alpha t}-1\g).
\end{align}

{\bf Proof of Theorem \ref{Qth1}.}\;Letting $\kappa\f(\frac{11}{8}2^{n}\g)^\alpha t_n=1$ in \eqref{hhy}, one has
\bbal
\De_n\f[\theta(t_n,\mathbf{x})-\theta_{0}(\mathbf{x})\g]&=(1-e^{-1})2^{-ns}\cos \f(\frac{11}{8}2^{n}x_1\g).
\end{align*}
Combining the above, we have
\bbal
\|\theta(t_n,\mathbf{x})-\theta_{0}(\mathbf{x})\|_{B^s_{p,\infty}(\T^2)}&\geq2^{ns}\f\|\De_n\f[\theta(t_n,\mathbf{x})-\theta_{0}(\mathbf{x})\g]\g\|_{L^p(\T^2)}
\\&=\f(1-e^{-1}\g)
\cdot\f\|\cos \f(\frac{11}{8}2^{n}x_1\g)\g\|_{L^p(\T^2)}\\
&=(2\pi)^{\frac{1}{p}} c_0\f(1-e^{-1}\g).
\end{align*}
Taking the $\limsup$, we deduce that
$$
\limsup_{t_n\to 0^+}\left\|\theta(t_n,\mathbf{x})-\theta_0\right\|_{B^s_{p,\infty}}\geq (2\pi)^{\frac{1}{p}} c_0\f(1-e^{-1}\g).
$$
Thus we complete the proof of Theorem \ref{Qth1}.

{\bf Proof of Theorem \ref{Qth2}.}\; We should notice that, the system \eqref{QG} with $\kappa=0$ has a unique steady-state solution $\theta^{0}(t,\mathbf{x})=\theta_0(\mathbf{x})$. Let $\kappa_n=t_n=\f(\frac{11}{8}2^{n}\g)^{-\frac{\alpha}{2}}$ in \eqref{hhy}, which gives that
\bbal
\De_n\f[\theta^{\kappa_n}(t_n,\mathbf{x})-\theta^{0}(t_n,\mathbf{x})\g]&=(1-e^{-1})2^{-ns}\cos \f(\frac{11}{8}2^{n}x_1\g),
\end{align*}
and thus we deduce that
$$
\limsup_{\kappa_n\to 0^+}\left\|\theta^{\kappa_n}(t_n,\mathbf{x})-\theta^{0}(t_n,\mathbf{x})\right\|_{B^s_{p,\infty}}\geq (2\pi)^{\frac{1}{p}} c_0\f(1-e^{-1}\g).
$$
Notice that $t_n\to 0^+$ and $\kappa_n\to0^+$ as $n\to\infty$, we complete the proof of Theorem \ref{Qth2}.
\section{Proof of Theorem \ref{Qth3}}\label{sec8}

\begin{lemma}\label{lm1} Assume that $s>0$ and $1 \leq p \leq \infty$. We define an initial data $\theta_0$  as follows
\bbal
\theta_0(\mathbf{x}):=
\sum^{\infty}_{j=3}2^{-js}\sin\f(\frac{11}{8}2^jx_1\g)+\sin x_2, \quad  \mathbf{x}=(x_1,x_2)\in \mathbb{T}^2.
\end{align*}
Then there exists some sufficiently large $n\in \mathbb{Z}^+$ and some positive constants $C,c$ such that
\bbal
&\f\|\theta_0\g\|_{L^\infty(\mathbb{T}^2)}+\f\|u_0\g\|_{L^\infty(\mathbb{T}^2)}\leq C,\\
&\f\|\theta_0\g\|_{B^s_{p,\infty}(\mathbb{T}^2)}\leq C,\\
&2^{ns}\f\|\De_{n}(u_0\cdot\nabla \theta_0)\g\|_{L^p(\mathbb{T}^2)}\geq c2^n.
\end{align*}
\end{lemma}
\begin{proof} Notice that
\bbal
\Lambda^{-1}\cos(\lambda x_i)=\frac{1}{2} \sum_{k \in \mathbb{Z}^2} \frac{1}{|k|}\f(\mathbf{1}_{\vec{\lambda}}(k)+\mathbf{1}_{-\vec{\lambda}}(k)\g)e^{\mathrm{i} x\cdot k}=\frac{1}{\lambda}\cos(\lambda x_i),\quad  i=1,2,
\end{align*}
then
\bbal
&u_0^{(1)}(x):=-\Lambda^{-1}\pa_{x_2}\theta_0(\mathbf{x})=-\cos x_2,\\
&u_0^{(2)}(x):=\Lambda^{-1}\pa_{x_1}\theta_0(\mathbf{x})=\sum^{\infty}_{j=3}2^{-js}\cos\f(\frac{11}{8}2^jx_1\g),
\end{align*}
which implies the first inequality. The second is obvious. Also, we have
\bbal
u_0\cdot\nabla \theta_0=
\cos x_2\cdot\sum^{\infty}_{j=3}2^{-js}\f(1-\frac{11}{8}2^j\g)\cos\f(\frac{11}{8}2^jx_1\g).
\end{align*}
From which, we have
\bbal
\De_{n}(u_0\cdot\nabla \theta_0)=\cos x_2\cdot2^{-ns}\f(1-\frac{11}{8}2^n\g)\cos\f(\frac{11}{8}2^nx_1\g),
\end{align*}
which in turn gives
\bbal
2^{ns}\f\|\De_{n}(u_0\cdot\nabla \theta_0)\g\|_{L^p(\mathbb{T}^2)}&= \f(\frac{11}{8}2^n-1\g)\left\|\cos x_2\right\|_{L^p([0,2\pi])}\left\|\cos\left(\frac{11}{8}2^{n}x_1\right)\right\|_{L^p([0,2\pi])}\\
&=\f(\frac{11}{8}2^n-1\g)\left(2\int^\pi_0|\cos x|^p\dd x\right)^{2/{p}}\geq c2^n,
\end{align*}
thus we complete the proof of Lemma \ref{lm1}.
\end{proof}
The following Lemma involves the nonlinear perturbation for the solution to the inviscid \eqref{QG}.
\begin{lemma}\label{pr3}
Assume that $s>0$ and $1 \leq p \leq \infty$. Let $(\theta,u)$ be a weak solution to the inviscid \eqref{QG} with initial data $\theta_0$. There exists a constant $C>0$ such that
\bbal
\left\|u\theta-u_0\theta_0\right\|_{\dot{B}^{s}_{p,\infty}}\leq C\f[\|\theta-\theta_0\|_{B^{s}_{p,\infty}}\|\theta_0\|_{L_{\mathcal{R}}^\infty}
+\|\theta-\theta_0\|_{L_{\mathcal{R}}^\infty}\f(\|\theta-\theta_0\|_{B^{s}_{p,\infty}}+\|\theta_0\|_{B^{s}_{p,\infty}}\g)\g].
\end{align*}
\end{lemma}

\begin{proof} We sketch the proof here. Due to the simple fact $u\theta-u_0\theta_0=(u-u_0)\theta+u_0(\theta-\theta_0)$, by the product law, one has
\bbal
\left\|u\theta-u_0\theta_0\right\|_{\dot{B}^{s}_{p,\infty}}&\leq C\Big(\|u-u_0\|_{\dot{B}^{s}_{p,\infty}}\|\theta\|_{L^\infty}+\|u-u_0\|_{L^\infty}\|\theta\|_{\dot{B}^{s}_{p,\infty}}\\
&\quad+\|\theta-\theta_0\|_{\dot{B}^{s}_{p,\infty}}\|u\|_{L^\infty}+\|\theta-\theta_0\|_{L^\infty}\|u\|_{\dot{B}^{s}_{p,\infty}}\Big)\\
&\leq C\f(\|\theta-\theta_0\|_{B^{s}_{p,\infty}}\|\theta\|_{L_{\mathcal{R}}^\infty}
+\|\theta-\theta_0\|_{L_{\mathcal{R}}^\infty}\|\theta\|_{B^{s}_{p,\infty}}\g),
\end{align*}
which implies the desired result.
\end{proof}

{\it Proof of Theorem \ref{Qth3}.}\;From \eqref{QG} and the Newton-Leibniz formula, it follows that
\bbal
\theta(t)-\theta_0&=-tu_0\cdot\nabla \theta_0-\int_0^t\div\left(u\theta-u_0\theta_0\right)\dd\tau.
\end{align*}
Using the triangle inequality, we have for some $C_*>0$
\bbal
\f\|\theta(t)-\theta_0\g\|_{B^{s}_{p,\infty}}
&\geq2^{{ns}}
\f\|\De_{n}\f(\theta(t)-\theta_0\g)\g\|_{L^p}\nonumber\\
&\geq t2^{{n}s}\f\|\De_{n}\big(u_0\cdot\nabla \theta_0\big)\g\|_{L^p}-C2^n\int_0^t2^{{n}s}\f\|\De_{n}(u\theta(\tau)-u_0\theta_0)\g\|_{L^p}\dd\tau\nonumber\\
&\geq 2^{n}\f[ct-C_*\int_0^t\f(\|\theta(\tau)-\theta_0\|_{B^{s}_{p,\infty}\cap L_{\mathcal{R}}^\infty}
+\|\theta(\tau)-\theta_0\|_{L_{\mathcal{R}}^\infty}\|\theta(\tau)-\theta_0\|_{B^{s}_{p,\infty}}\g)\dd\tau\g],
\end{align*}
where we have used Lemmas \ref{lm1}-\ref{pr3}. In fact, we have for some $C_*>0$
\bbal
\left\|u\theta-u_0\theta_0\right\|_{\dot{B}^{s}_{p,\infty}}&\leq C_*\f(\|\theta-\theta_0\|_{B^{s}_{p,\infty}\cap L_{\mathcal{R}}^\infty}
+\|\theta-\theta_0\|_{L_{\mathcal{R}}^\infty}\|\theta-\theta_0\|_{B^{s}_{p,\infty}}\g).
\end{align*}

Suppose that \eqref{imp2} were not true for some $s>0$ and $1 \leq p \leq \infty$. We pick $t_0\in(0, t]$ such that for all $\tau\in[0,t_0]$
$$
\f\|\theta(\tau)-\theta_0\g\|_{B^{s}_{p,\infty}\cap L_{\mathcal{R}}^{\infty}}<\min \left\{1, \frac{c}{4C_*}\right\}.
$$
Thus, picking $\tau_n=M2^{-n}$ with large $M>0$, we have
\bbal
1>\f\|\theta(\tau_n)-\theta_0\g\|_{B^{s}_{p,\infty}}\geq  \frac{c}{2}M,
\end{align*}
which leads to a contraction for large $M$. This completes the proof of Theorem \ref{Qth3}.

\section*{Declarations}
\noindent\textbf{Data Availability}\;
 No data was used for the research described in the article.

\vspace*{1em}
\noindent\textbf{Conflict of interest}\;
The authors declare that they have no conflict of interest.
\vspace*{1em}

\noindent\textbf{Funding}\;
J. Li was supported by the Ganpo Talents Project of Jiangxi Province (gpyc20240069) and Natural Science Foundation of Jiangxi Province (20252BAC210004 and 20252BAC240186). X. Wu was supported by the Natural Science Foundation of Henan Province (252300421774).
\section*{Acknowledgements}
The authors would like to thank the anonymous referees for valuable comments and suggestions which greatly improved the presentation of this paper.


\begin{thebibliography}{99}
\linespread{0}\addtolength{\itemsep}{-1.0ex}

\bibitem{BCD} H. Bahouri, J.-Y. Chemin, R. Danchin, Fourier Analysis and Nonlinear Partial Differential Equations, Grundlehren der Mathematischen Wissenschaften, 343, Springer, New York (2011).
\bibitem{Be} F. Bernicot, T. Elgindi, S. Keraani, On the inviscid limit of the 2D Navier-Stokes equations with
vorticity belonging to BMO-type spaces, Ann. I. H. Poincar\'{e}-AN., 33 (2016), 597-619.
\bibitem{B1} J. Bourgain, D. Li, Strong ill-posedness of the incompressible Euler equation in borderline Sobolev spaces, Invent. Math., 201 (2015), 97-157.
 \bibitem{B2} J. Bourgain, D. Li, Strong ill-posedness of the incompressible Euler equation in integer $C^m$ spaces, Geom. Funct. Anal., 25 (2015), 1-86.
 \bibitem{chem} J.-Y. Chemin, A remark on the inviscid limit for two-dimensional incompressible fluids, Commun. Partial Differ. Equ., 21 (11-12) (1996), 1771-1779.



\bibitem{Chae1} D. Chae, On the well-posedness of the Euler equations in the Triebel-Lizorkin spaces, Commun. Pure Appl. Math., 55 (2002), 654-678.
\bibitem{Chae2} D. Chae, On the Euler equations in the critical Triebel-Lizorkin spaces, Arch. Ration. Mech. Anal., 170 (2003), 185-210.
\bibitem{Chae3} D. Chae, Local existence and blowup criterion for the Euler equations in the Besov spaces, Asymptot. Anal., 38 (2004), 339-358.
\bibitem{CMZ} Q. Chen, C. Miao, Z. Zhang, On the well-posedness of the ideal MHD equations in the Triebel-Lizorkin spaces,
Arch. Ration. Mech. Anal., 195 (2010), 561-578.
\bibitem{CS} A. Cheskidov, R. Shvydkoy, Ill-posedness of the basic equations of fluid dynamics in Besov spaces, Proc. Amer. Math. Soc., 138 (2010), 1059-1067.
\bibitem{Ci} G. Ciampa, G. Crippa, S. Spirito, Strong convergence of the vorticity for the 2D Euler equations in the inviscid limit, Arch. Rational Mech. Anal., 240 (2021), 295-326.

\bibitem{P2} P. Constantin, T. Drivas, T. Elgindi, Inviscid limit of vorticity distributions in Yudovich class, Commun. Pure Appl. Math., 75 (2022), 0060-0082.
\bibitem{C17} P. Constantin, T. Elgindi, M. Ignatova, V. Vicol, Remarks on the inviscid limit for the Navier-Stokes
equations for uniformly bounded velocity fields, SIAM J. Math. Anal., 49 (2017), 1932-1946.
\bibitem{CMT94} P. Constantin, A. Majda, E. Tabak, Formation of strong fronts in the 2D quasi-geostrophic thermal active scalar, Nonlinearity, 7 (1994), 1495-1533.
\bibitem{CKV}  P. Constantin, I. Kukavica, V. Vicol, On the inviscid limit of the Navier-Stokes equations, Proc. Amer. Math. Soc., 143 (2015), 3075-3090.

\bibitem{CV}  P. Constantin, V. Vicol, Remarks on high Reynolds numbers hydrodynamics and the inviscid limit, J. Nonlinear Sci., 28 (2018), 711-724.
\bibitem{CM22ADV} D. C\'{o}rdoba, L. Mart\'{\i}nez-Zoroa, Non existence and strong ill-posedness in $C^k$ and Sobolev spaces for SQG, Adv. Math., 407 (2022), 108570.
\bibitem{CM24CMP} D. C\'{o}rdoba, L. Mart\'{\i}nez-Zoroa, Non-existence and strong ill-posedness in $C^{k,\,\beta}$ for the generalized surface quasigeostrophic equation, Commun. Math. Phys., 405 (2024), 170.
\bibitem{Dai20} Y. Dai, W. Hu, J. Wu, B. Xiao, The Littlewood-Paley decomposition for periodic functions and applications to the Boussinesq equations, Anal. Appl.(Singap.), 18 (2020), 639-682.
\bibitem{Elgindi20arma} T. M. Elgindi, N. Masmoudi, $L^\infty$ ill-posedness for a class of equations arising in hydrodynamics, Arch. Ration. Mech. Anal. 235:3 (2020), 1979-2025.
\bibitem{JK24APDE} I.-J. Jeong, J. Kim, Strong ill-posedness for SQG in critical Sobolev spaces, Anal. PDE., 17 (2024), 133-170.
\bibitem{guo} Z. Guo, J. Li, Z. Yin, Local well-posedness of the incompressible Euler equations in $B^1_{\infty,1}$ and the inviscid limit of the Navier-Stokes equations, J. Funct. Anal., 276 (2019), 2821-2830.
\bibitem{HK} T. Hmidi, S. Kerrani, Inviscid limit for the two-dimensional N-S system in a critical Besov space, Asymptot. Anals., 53 (2007), 125-138.
\bibitem{HK1} T. Hmidi, S. Kerrani, Incompressible viscous flows in borderline Besov spaces, Arch. Rational Mech. Anal., 189 (2008), 283-300.


\bibitem{Kato} T. Kato, Nonstationary flows of viscous and ideal fluids in $\mathbb{R}^3$, J. Funct. Anal., 9 (1972), 296-305.

\bibitem{KatoP} T. Kato, G. Ponce, Commutator estimates and the Euler and Navier-Stokes equations, Commun. Pure Appl. Math., 41 (1988), 891-907.
\bibitem{KNV} A. Kiselev, F. Nazarov, A. Volberg, Global well-posedness for the critical 2D dissipative quasi-geostrophic equation, Invent. math., 167 (2007), 445-453.
\bibitem{AA2} M. Kohr, R. Precup, Localization of energies in Navier-Stokes models with reaction terms, Anal. Appl.(Singap.), 22 (2024), 1053-1073.

\bibitem{Lyz1} J. Li, Y. Yu, W. Zhu, Ill-posedness for the Euler equations in Besov spaces, Nonlinear Anal. Real World Appl., 74 (2023) 103941.
 \bibitem{BULL} J. Li, Y. Yu, N. Zhu, On the continuous properties for the 3D incompressible rotating Euler equations, Bull. Sci. Math., 211 (2026), Paper No. 103821.
\bibitem{Majda} A. Majda, Vorticity and the mathematical theory of an incompressible fluid flow, Commun. Pure Appl. Math., 39 (1986), 187-220.

\bibitem{M} N. Masmoudi, Remarks about the inviscid limit of the Navier-Stokes system, Commun. Math. Phys., 270 (2007), 777-788.
\bibitem{Sawano} M. Masty{\l}o, Y. Sawano, The Navier-Stokes equation in some general Banach spaces, Anal. Appl.(Singap.), 24 (2026), 1055-1121.
\bibitem{MY} G. Misio{\l}ek, T. Yoneda, Ill-posedness examples for the quasi-geostrophic and the Euler equations, Contemp. Math., 584
 Amer. Math. Soc., (2012), 251-258.
\bibitem{MY2} G. Misio{\l}ek, T. Yoneda, Continuity of the solution map of the Euler equations in H\"{o}lder spaces and weak norm
inflation in Besov spaces, Trans. Amer. Math. Soc., 370 (2018), 4709-4730.

\bibitem{Pak} H.C. Pak, Y.J. Park, Existence of solutions for the Euler equations in a critical Besov space $B^1_{\infty,1}(\R^n)$, Commun.
Partial Differ. Equ., 29(7-8) (2004), 1149-1166.

\bibitem{Swann} H.S.G. Swann, The convergence with vanishing viscosity of nonstationary Navier-Stokes flow to ideal flow in $\mathbb{R}^3$, Trans. Amer. Math. Soc., 157 (1971), 373-397.
\bibitem{YuL} Y. Yu, J. Li, Non-convergence of the Navier-Stokes equations toward the Euler equations in the endpoint Besov spaces, Appl. Math. Optim., 92 (2025), Paper No.32.

\bibitem{YuLf} Y. Yu, F. Liu, Ill-posedness for the Navier-Stokes and Euler equations in Besov spaces, Appl. Math., 69 (2024), 757-767.

\end{thebibliography}
\end{document}